\documentclass[preprint,3p, times, 12pt]{elsarticle}

\usepackage{amssymb,amsmath,mathrsfs,framed}
\usepackage{hyperref,color}
\usepackage{paralist}
\usepackage[all]{xy}

\journal{}
\newtheorem{thm}{Theorem}[section]

\newtheorem{example}{Example}[section]

\newproof{pf}{Proof}
\newdefinition{rmk}{Remark}
\newcommand{\bF}{\mathbb{F}}
\newcommand{\Finf}{\mathbb{F}_q^{\cup\{\infty\}}}
\newcommand{\zz}{z}
\newcommand{\PGL}{\mathrm{PGL}}
\newcommand{\PG}{\mathrm{PG}}
\newcommand{\PSL}{\mathrm{PSL}}
\newcommand{\SL}{\mathrm{SL}}
\newcommand{\PML}{\mathrm{P}\Gamma\mathrm{L}}
\newcommand{\ML}{\Gamma\mathrm{L}}
\newcommand{\GL}{\mathrm{GL}}
\newcommand{\POm}{\mathrm{P}\Omega}

\newcommand{\fX}{\mathfrak{X}}
\newcommand{\SO}{\mathrm{SO}}
\newcommand{\GO}{\mathrm{GO}}
\newcommand{\FT}{\mathrm{FT}}
\newcommand{\Sym}{\mathrm{Sym}}
\newcommand{\Mq}{\mathrm{M}}
\newcommand{\col}[1]{\left(\begin{array}{c} #1 \end{array}\right)}

\begin{document}

\begin{frontmatter}

\title{Fissioned triangular schemes via sharply 3-transitive groups}

\author[JM]{Jianmin Ma } 
\ead{Jianmin.Ma@gmail.com}
\address[JM]{College of Math \&  Info. Science, Hebei Normal University, Shijiazhuang, 050016, China}
\author[KW]{Kaishun Wang}
\ead{wangks@bnu.edu.cn}
\address[KW]{Sch. Math. Sci. \& Lab. Math. Com. Sys., Beijing Normal University, Beijing  100875, China}


\begin{abstract}
In  [D. de Caen, E.R. van Dam. Fissioned triangular schemes via the cross-ratio,  { Europ. J. Combin.}, 22 (2001) 297-301], de Caen and  van Dam constructed a fission scheme $\FT(q+1)$ of the triangular scheme on $\PG(1,q)$. This fission scheme comes from the naturally induced action of $\PGL(2,q)$ on the 2-element subsets of $\PG(1,q)$. The group $\PGL(2,q)$  is one of two infinite families of finite sharply 3-transitive groups.  The other such family $\Mq(q)$ is a ``twisted'' version of $\PGL(2,q)$, where $q$ is an even power of an odd prime. The group $\PSL(2,q)$ is the intersection of $\PGL(2,q)$ and $\Mq(q)$.  
In this paper, we investigate the association schemes coming from the actions of $\PSL(2,q)$, $\Mq(q)$ and $\PML(2,q)$, respectively. Through the conic model introduced  in [H.D.L. Hollmann,  Q. Xiang.  Association schemes  from the actions of
$\PGL(2, q ) $ fixing a nonsingular conic,  { J. Algebraic Combin.}, 24 (2006) 157-193],  we introduce an 
 embedding of  $\PML(2,q)$ into $\PML(3,q)$. For each of  the three groups mentioned above, this embedding produces  two more isomorphic association schemes:  one on  hyperbolic lines and the other on hyperbolic points (via a null parity)  in a 3-dimensional orthogonal geometry.    This embedding enables us to   treat these  three isomorphic association schemes simultaneously.

\end{abstract}


\end{frontmatter}
\section{Introduction} \label{s:introduction}
Let $X$ be a finite set with cardinality $n \ge 2$  and $\mathbf{R}=\{R_0, R_1, \ldots , R_d\}$ be a set 
of binary relations on $X$.  ${\mathfrak X} = (X,\mathbf{R})$ is called an \emph{association scheme with $d$ classes} 
(a {\em  $d$-class association scheme}, or simply, a \emph{scheme}) if the following axioms are satisfied:
\renewcommand{\theenumi}{\roman{enumi}}
\renewcommand{\labelenumi}{(\theenumi)}
\begin{enumerate}
\item 
$\mathbf{R}$ is a partition of $X\times X$ and $R_0 = \{(x,x)\;|\;x \in X\}$ is the diagonal relation;
\item  
For $i=0,1,\dots, d$, the inverse ${^tR_i} = \{(y, x) \vert (x, y) \in  R_i\}$ of $R_i$ 
is also among the relations:  ${^tR_i} = R_{i'}$ for some $i'$ ($0 \le  i' \le d$);
\item 
For any triple of $i, j, k = 0, 1, \dots, d$, there exists an integer $p^k_{ij}$ such that for all
$(x,y) \in R_k$,
\[
 |\{z \in  X\;|\; (x, z) \in  R_i, (z, y)\in  R_j\}| = p^k_{ij}.
\] 
\end{enumerate}  
The integers $p^k_{ij}$ are called the \emph{intersection numbers}. 
The integer $k_i = p^0_{i i'} $ is called the \emph{valency} of $R_i$. 
In fact, for any $x\in X$, $k_i = |\{y\in X \vert\; (x,y) \in R_i\}|$.

If  $ p^k_{ij} =p^k_{j\,i}$ for all $i, j, k$, ${\mathfrak X}$ is called \emph{commutative}.
 A relation $R_i$ is said to be \emph{symmetric} if  $R_i = {^tR}_i$.  The scheme $\mathfrak{X}$ is called \emph{symmetric} if all relations $R_i$ are symmetric.  
A partition 
$\Lambda_0, \Lambda_1,\ldots, \Lambda_e$ of the index set $\{0,1,\dots,d\} $
is called \emph{admissible} \cite{ItM91} if $\Lambda_0 =\{0\}$, $\Lambda_i \ne \emptyset$ and 
${\Lambda'_i} = \Lambda_j$  for some $j\ (1\le i, j\le e)$, where  
${\Lambda'} =\{\alpha' | \alpha \in \Lambda\}$. Let
$R_{\Lambda_i} = \cup_{\alpha\in \Lambda_i} R_\alpha$.  If  
$\mathfrak{Y}= (X, \{R_{\Lambda_i}\}_{i=0}^e)$ is an association scheme, it is called
 a \emph{fusion scheme} of ${\mathfrak X}$, and $\mathfrak{X}$ is called a {\em fission scheme} of $\mathfrak{Y}$. 

We mention  a typical example of association scheme (see \cite[Section 2.1]{BI} for details). If group $G$ acts transitively  on a finite set $X$, then the orbits of  the induced diagonal action of $G$ on $X\times X$ form an association scheme, denoted by $\mathfrak{X}(G, X)$. For any $x\in X$, the orbits of $G_x$ on $X$ are in one-to-one correspondence with those of $G$ on $X\times X$, where $G_x$ is the stabilizer of $x$ in $G$.    The scheme $\mathfrak{X}(G, X)$ is commutative if and  only if the permutation representation of $G$ on $X$ is {\em multiplicity-free}, and it is symmetric if and only if $G$ acts {\em generously transitively} on $X$.

By a theorem of Zassenhaus \cite[Section 7.1]{Robinson}, there are two infinite families of finite sharply 3-transitive groups, and both are subgroups of the projective semilinear group $\PML(2,q)$. If $G$ is a sharply 3-transitive group of permutations on a finite set $X$, then $X$ can be identified with the projective line $\PG(1,q)$ for some prime power $q$ and $G$ is one of the following: 
\begin{enumerate}
\item [(1)]
$G$ is  the projective general linear group $\PGL(2,q)$ in its natural action on $\PG(1,q)$;
\item  [(2)]
$q=p^{2f}$ for some odd prime $p$ and a positive integer $f$, and if $\sigma$ is the unique involution in $\mathrm{Aut} (\bF_q)$, then $G$ is the group $\Mq(q) = S\cup T$, where  
\begin{align*}
S  & =  \left \{\lambda \mapsto \frac{a\lambda+b}{c\lambda +d},  ~\right\vert 
                 \left. ad-bc \mbox{ is a square in } \bF_q^* \right\},
\\                 
T & = \left \{\lambda \mapsto  \frac{a\lambda^\sigma +b}{c\lambda^\sigma +d} ~\right\vert
 \left. ad-bc \mbox{ is a non-square  in }  \bF_q^* 	\right\}.
\end{align*}
\end{enumerate}

In this paper,  we shall investigate the association schemes determined by the following  subgroups of the symmetric group
$\mbox{Sym}(q+1)$:
$$
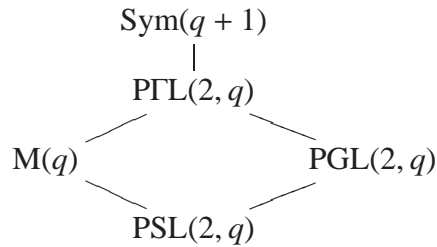
\xymatrix@-1.2pc{
& \Sym(q +1) & \\
& \PML(2, q) \ar@{-}[u]& \\
\Mq(q) \ar@{-}[ur] & & \PGL(2,q) \ar@{-}[ul] \\
& \PSL(2,q) \ar@{-}[ul] \ar@{-}[ur] &
}
$$
These groups  are permutation groups on $\PG(1,q)$ and they acts transitively on $\Omega$,  the collection of 2-element subsets of  $\PG(1,q)$.  Hence, each of these groups determines an association scheme on $\Omega$:    $\mathrm{Sym}(q+1)$  determines the triangular scheme $\mathrm{T}(q+1)$, and $\PGL(2,q)$ determines a symmetric fission scheme of   
$\mathrm{T}(q+1)$.  De Caen and van Dam  \cite{DeCaen99} described $\fX(\PGL(2,q), \Omega)$ via the cross-ratio.  See Subsection \ref{s:FT}. The association schemes from the actions of the remaining three groups are the main object of the present paper. 

The isomorphism $\PGL(2,q)\simeq \SO(3,q)$ introduces a conic model (see Subsection \ref{s:O3geom}). The action $\PGL(2,q)$ on $\Omega$ is equivalent to  that of $\SO(3,q)$ on the set of hyperbolic lines $\mathcal{L}_+$ in a 3-dimensional orthogonal geometry.  This model has been used heavily in  \cite{Ebert01,Ebert02,Xiang03}.  This conic model allows us to use an embedding $\rho$ of $\PML(2,q)$ into $\PML(3,q)$ and the image of this embedding fixes this conic as a set. As a result,   the image of $\PSL(2,q)$ acts transitively on the hyperbolic lines $\mathcal{L}_+$ and so do the images of $\Mq(q)$ and $\PML(2,q)$.   

The conic above also introduces a (null) polarity $\perp$ such that  the action of any  subgroup of $\PML(2,q)$ on the hyperbolic lines $\mathcal{L}_+$ is equivalent to that of its embedded image on  the hyperbolic points $\mathcal{L}_+^\perp$;  see Subsection \ref{s:q-odd}. 
Therefore, we do not  need to distinguish actions on hyperbolic lines and points in certain calculation. 
 The advantage of this conic model allows us to treat (isomorphic) association schemes on $\Omega$,  $\mathcal{L}_+$ and 
 $\mathcal{L}_+^\perp$ all at once.

This paper is organized as follows.  In Section \ref{s:isomorphism}, we introduce the conic model, the polarity and the embedding of $\PML(2,q)$ into $\PML(3,q)$ mentioned above, and establish a few results on transitivity. In Section 
\ref{s:lattice}, we determine three fission schemes of the triangular scheme $\mathrm{T}(q+1)$ via their isomorphic association schemes.  

 We refer to \cite{BI} for undefined terms   and the basic theory of  association schemes and to  \cite{Cameron00, Taylor92, Wan02} for missing definitions and notation about various groups in this paper.

\section{The 3-dimensional orthogonal geometry and an embedding of $\PML(2,q)$} \label{s:isomorphism}
In this section, we introduce the 3-dimensional orthogonal geometry, a conic model, and various groups related to this conic.  

\subsection{The 3-dimensional orthogonal geometry} \label{s:O3geom}

Let $\bF_q$ be a finite field with $q$ elements. Let $V = \bF_q^3$ be a 3-dimensional vector space  over $\bF_q$ equipped with a non-degenerate quadratic form $Q$.  {\em The general orthogonal group} $\GO(V)$ is the isometry group of $V$ with respect to $Q$:
$$
\GO(V) = \{ A \in \GL(V) ~\vert~ Q(A(x_0, x_1,x_2)^T) = Q(x_0, x_1, x_2)\}.
$$
{\em The special linear orthogonal group} $\SO(V)$ is the intersection of $\GO(V)$ and $\SL(V)$. 
It is  also standard to write $\GO(3,q)$ for  $\GO(V)$, and    $\SO(3,q)$ for  $\SO(V)$ when the underlying field is $\bF_q$. 

The projective plane $\PG(2,q)$ on $V$ has as  points  the 1-dimensional subspaces of $V$ and as lines (hyperplanes) the 2-dimensional subspaces.  Any point $P$ is spanned by a nonzero vector $v= (v_0, v_1, v_2)$. 
Another vector $u = (u_0, u_1, u_2)$ spans $P$ if and only if $v = \xi u$ for some $\xi\in \bF_q^*$. 
So we use  $(v_0 :  v_1 : v_2)$ to denote the point $P$. 
A point  $(v_0 :  v_1 : v_2)$ is called \emph{singular} if $Q(v_0, v_1, v_2) =0$.

In the rest of this paper,  we fix the quadratic form $Q(x_0, x_1, x_2) = x_1^2 - x_0x_2$. Let 
\begin{equation} \label{e:conic}
\mathscr{O} = \{(\xi^2: \xi : 1) ~\vert~ \xi \in \bF_q\} \cup \{(1:0:0)\}.
\end{equation}
    Then $\mathscr{O}$ is a \emph{conic}. 
 Let $ \infty = (1:0:0).$ We use $P_\xi$ and $P_\infty$ to denote the points  $(\xi^2: \xi : 1)$ and $(1:0:0)$, respectively.
 No three distinct points of $\mathscr{O}$ can be on a line. 
 Therefore, any line $\ell$  intersects $\mathscr{O}$ at most 2 points. Accordingly,   $\ell$ is called 
{\em hyperbolic} or {\em secant}  if $|\ell \cap \mathscr{O}| = 2$, {\em  tangent}  if $|\ell \cap \mathscr{O}| = 1$, or {\em elliptic} or {\em exterior} if $|\ell \cap \mathscr{O}| = 0$.   We denote by $\mathcal{L}_+$, $\mathcal{L}_{0}$ and $\mathcal{L}_-$  the set of all hyperbolic,  tangent and elliptic lines, respectively.

The following result is well known (e.g. see \cite[Theorem 11.6]{Taylor92}), from which Theorem  \ref{t:SO3} follows immediately. 

\begin{thm}  \label{t:O3}
If $V$ is an orthogonal geometry of dimension 3 and Witt index 1 over $\bF_q$,  then 
$\GO(V) \simeq \{\pm 1\} \times \SO(V)$, $\SO(V) \simeq \PGL(2,q),$ and $\SO(V)$ acts triply transitively on the  set of all singular points.
\end{thm}

\begin{thm}\label{t:SO3}
The group
$\SO(V)$ acts  generously transitively on $\mathcal{L}_{\epsilon}$ for $\epsilon \in \{+, 0\}$.  
\end{thm}

We remark that $\SO(V)$ also acts  generously transitively on $\mathcal{L}_{-}$;  see  \cite{Shen84}.
Therefore,  the action $(\SO(V),  \mathcal{L}_{\epsilon})$ determines a symmetric association scheme, denoted by $\mathfrak{X}(\SO(3,q),  \mathcal{L}_{\epsilon}), \epsilon=0,\pm$.

\subsection{An embedding of $\PML(2,q)$ into $\PML(3,q)$} \label{s:embedding}
In this subsection, we describe an embedding of   $\PML(2,q)$ into $\PML(3,q)$, which gives 
an isomorphism  $\PGL(2,q)\simeq \SO(3,q)$ mentioned in Theorem \ref{t:O3} (cf. \cite[Section 6.1]{Cameron00} and 
\cite[Section 3]{Hollmann}).   While this is a folklore, a detailed account is provided here to prepare for later sections.   

We first describe the action of $\PML(2,q)$  on the projective line $\PG(1,q)$.  Let 
$$
\PG(1,q) = \{ (\xi: 1) ~\vert~ \xi \in \bF_q\} \cup \left\{ (1:0) \right\} 
$$
and   $\Finf = \bF_q \cup \{\infty\}$, the union of $\bF_q$ and $\{\infty\}$.
Now we may identify  $\PG(1,q)$ with $\Finf$ by the map: $(\xi:1) \leftrightarrow \xi$, 
$(1:0) \leftrightarrow \infty$.  The group $\ML(2,q)$  is the wreath product of $\GL(2,q)$ and $\mathrm{Aut}( \bF_q)$, and
  $\PML(2,q)$ is  the group induced on the the projective line $\PG(1,q)$ by  $\ML(2,q)$. 
Each element of $\PML(2,q)$ is induced by  a pair $(A, \tau)$ with $A= \begin{pmatrix} a & b \\ c & d \end{pmatrix} \in \GL(2,q)$ and 
$\tau\in \mathrm{Aut}(\bF_q)$, which acts on $\PG(1,q)$ as follows:
$$
  \xi \mapsto A(\xi^\tau) :=  \frac{a \xi^\tau +b}{c\xi^{\tau} +d}  \quad \mbox{ for all } \xi\in \Finf .
$$
The expressions involving $\infty$ are evaluated by  standard limit rules: e.g., $\infty^\tau = \infty$, $A(\infty) = a/c$ if $c\ne 0$ and $A(\infty) = \infty$ if $c=0$.

Consider the vector space $W$ of quadratic forms on $\bF_q^2$: $q(x,y)= u x^2 + v xy + w y^2$. So, $\mbox{dim} W =3$. 
 The group $\ML(2,q)$ acts on $W$ as follows: for each $(A, \tau)$ in $\ML(2,q)$, 
\begin{equation}\label{pgl2qactQ}
q(x,y)   \mapsto q(A(x^\tau,y^\tau)^T)  \mbox{~for all $q \in W$}. 
\end{equation}
For brevity, we write $q(A(x^\tau,y^\tau)) = q(A(x^\tau,y^\tau)^T) $. 
If  $ A = \begin{pmatrix} a & b \\ c & d \end{pmatrix} \in \GL(2,q)$,  then 
\begin{equation*}
\begin{split}
q(A(x^\tau,y^\tau))  = & 
u (ax^\tau +by^\tau)^2 + v (ax^\tau + by^\tau) (cx^\tau + dy^\tau) + w (cx^\tau +dy^\tau)^2  \\
                     = & x^{2\tau} (u a^2 + v ac + w c^2) + x^\tau y^\tau (2uab + v(ad+bc) + 2w cd) + y^{2\tau} ( u b^2 + v bd + w d^2) \\
                      = & 
                     \begin{pmatrix} x^{2\tau}  & x^\tau y^\tau & y^{2\tau} \end{pmatrix} 
                       \begin{pmatrix}  a^2 & 2ab & b^2 \\  ac & ad +bc & bd \\ c^2 & 2cd & d^2 \end{pmatrix}
                       \begin{pmatrix} u \\  v \\ w \end{pmatrix}. 
\end{split}
\end{equation*}
Let 
\begin{equation} \label{e:so3q}
\rho(A) = \begin{pmatrix}  a^2 & 2ab & b^2 \\  ac & ad +bc & bd \\ c^2 & 2cd & d^2 \end{pmatrix}. 
\end{equation}
Since $\det \rho(A) = \det(A)^3$, it can be verified that  (\ref{pgl2qactQ}) induces a homomorphism from $\ML(2,q)$ to $\ML(3,q)$.     In particular, if $\tau$ is the identity automorphism,  $\rho$ is a homomorphism of $\GL(2,q)$ into $\GL(3,q)$. The kernel of $\rho$ is $\{\pm I\}$, where $I$ is the $2\times 2$ identity matrix . Let $Z =  \{ cI | c\in \bF_q^*\}$. Since $\PGL(2,q) = \GL(2,q)/Z$, $\rho$ embeds $\PGL(2,q)$ into $\PGL(3,q)$.  
{By abuse of notation, we also use $\rho$ to denote the homomorphism induced by
 (\ref{pgl2qactQ}).
 }
   Thus, $\rho$ embeds $\PML(2,q)$ into $\PML(3,q)$.

Since 
\begin{gather}
 \begin{pmatrix}  a^2 & 2ab & b^2 \\  ac & ad +bc & bd \\ c^2 & 2cd & d^2 \end{pmatrix}
 \begin{pmatrix} \xi^{2\tau} \\ \xi^\tau\\ 1 \end{pmatrix}
= \begin{pmatrix} (a\xi^\tau +b)^2 \\ (a\xi^\tau+b)(c\xi^\tau + d) \\ (c\xi^\tau +d )^2 \end{pmatrix},
\end{gather}
$(\rho(A),\tau)$ maps $(\xi^2: \xi:1)$ to 
$\left(\left(\frac{a\xi^\tau+b}{c\xi^\tau + d}\right)^2,  \frac{a\xi^\tau+b}{c\xi^\tau + d}:1\right)$ if 
$c\xi^\tau +d \ne 0$ and  to $(1:0:0)$ if $c\xi^\tau +d =0$. Similarly, $(\rho(A),\tau)$ maps  $(1:0:0)$ to $(1:0:0)$ if $c=0,$
 and to  $((a/c)^2, {a/c}, 1)$  if $c\ne 0$. 
Therefore,  $\rho(\PML(2,q))$ fixes the set $\mathscr{O}$.

Now we equip the space $W$ with the quadratic form $Q(u,v, w) = v^2 - u w$. By (\ref{pgl2qactQ}), $\GL(2,q)$ acts on $W$.
 Note that 
$$ Q(q(A(x,y)) = \det (A)^2 Q(q(x,y)),
$$
 i.e.,$Q$  is multiplied by  a factor $\mbox{det}(A)^2$ by the action of $\rho(A)$.
Hence, we find a subgroup of $\GO(3,q)$ which is isomorphic to $\SL^{\pm}(2,q) /\{\pm I\}$, where 
 $\SL^{\pm}(2,q)$ is the group of matrices with determinant $\pm 1$. 
In general,  $\PSL(2,q) \simeq \POm(3,q)$ and they are simple if $q>3$. As a matter of fact, 
$\rho(\PSL(2,q)) = \POm(3,q)$.

Note that $\SL^{\pm}(2,q) /\{\pm I\}$ and $\PGL(2,q)$ have the same size. For $q$  even,  $\SL^{\pm}(2,q) /\{\pm I\} = \SL(2,q)= \PGL(2,q)$.  For $q$ odd, $\SL^{\pm}(2,q) /\{\pm I\}$  and $ \PGL(2,q)$ are not isomorphic.  In this case, 
$\rho(\SL^{\pm}(2,q) /\{\pm I\})$ is a subgroup $\mathrm{PGO}(3,q)$ of index 2, but it is not isomorphic to $\mathrm{PSO}(3,q)$.

Define a map $f$ from $\PG(1,q)$ to $\mathscr{O}$ as follows:
\begin{equation*}
f(\infty) = \infty ~\mbox{and}~ 
f  (\xi, 1)  \mapsto (\xi^2 : \xi : 1)  \mbox{ for all $\xi \in \bF_q$}.
\end{equation*} 
We have 
\begin{equation} \label{e:quivact}
f(A(\xi)) = \rho(A) f(\xi)^T.
\end{equation}
The action of $\PGL(2,q)$ on $\PG(1,q)$ is equivalent to that of $\SO(3,q)$ on $\mathscr{O}$. 

\subsection{$\FT(q+1)$ and $\mathfrak{X}(\SO(3,q), \mathcal{L}_{\epsilon})_{\epsilon = \pm}$}
\label{s:FT}

The {\em triangular scheme} $\mathrm{T}(n)$ comes from the action of  $\mathrm{Sym}(n)$ on the collection of 2-element subsets of an $n$-set ($n\ge 4$): for any such subsets $x, y$, $(x,y) \in R_i$ if their intersection has size $2-i$.  $\mathrm{T}(n)$ is a   2-class Johnson scheme. 
It is well known that  $\PGL(2,q)$ acts  sharply 3-transitively on $\PG(1,q)$. So $\PGL(2,q)$ acts generously transitively on the collection of 2-element subsets of $\PG(1,q)$ and thus this action determines a symmetric association scheme. In the paper \cite{DeCaen99}, de Caen and van Dam described this scheme as  a fission scheme of the triangular scheme on $\PG(1,q)$ using the cross ratio, denoted by $\FT(q+1)$.

We  know from the previous subsection that the action $\PGL(2,q)$ on $\PG(1, q)$ is equivalent to that of 
$\SO(3,q)$ on $\mathscr{O}$. So the schemes   $\FT(q+1)$ is isomorphic to the scheme coming from the action of  
$\SO(3,q)$ on the 2-element subsets of $\mathscr{O}$.  By definition,  the hyperbolic lines $\mathcal{L}_+ $  are in one-to-one correspondence with the 2-element subsets of $\mathscr{O}$. Hence,  $\FT(q+1)$ is isomorphic to $\mathfrak{X}(\SO(3,q), \mathcal{L}_+)$.  In \cite{Hollmann},  Hollmann and Xiang investigated
$\mathfrak{X}(\SO(3,q), \mathcal{L}_{\epsilon})$ for $\epsilon =\pm$ using the cross ratio. In fact, they studied the coherent configuration from the action of $\PGL(2,q)$ $(\simeq \SO(3,q)$) on $\mathcal{L}_{+} \cup \mathcal{L}_{-}$.  For the case $q$ being a power of 4, they  calculated the intersection numbers  for this coherent configuration and constructed another coherent configuration on  $\mathcal{L}_{+} \cup \mathcal{L}_{-}$ by merging certain relations.  In this new coherent configuration,  the fiber $\mathcal{L}_+$ supports a 4-class scheme  and $\mathcal{L}_{-}$ supports  a 3-class scheme.

The 4-class  scheme on  $\mathcal{L}_+$  has an  interesting history. Its existence as a fusion scheme of   $\FT(q+1)$ was first conjectured by de Caen and van Dam \cite{DeCaen99}.  Tanaka \cite{Tanaka02} proved this conjecture  with character-theoretic  method.  Ebert et al.  proved it with a geometric interpretation in \cite{Ebert01} and with a direct calculation of intersection numbers  in \cite{Ebert02}.   Xiang gave a summary of these work in \cite{Xiang03}.

The  3-class scheme is also very interesting.  H. Tanaka observed that  this 3-class scheme on $\mathcal{L}_{-}$   gives a family of primitive schemes having the same parameters as  the first infinite  family of $Q$-polynomial but not $P$-polynomial association schemes by Penttila and Williford  \cite{Penttila10}. (The reference to Hollmann and Xiang  is incorrectly cited in  \cite{Penttila10}, and it should be \cite{Hollmann}. )

We conclude this subsection with a description of $\FT(q+1)$ following \cite{DeCaen99}.  It has as vertices all  2-element subsets of $\PG(1,q) $ as and its relations are  as follows:
\begin{equation}\label{e:rels}
\left.\begin{array}{cl}
R_0 & = \left\{ ( \{\xi,\gamma\}, \{\xi,\gamma\} ) ~\vert~  \xi, \gamma \in \Finf \right\}  
\\\noalign{\medskip}
R_1 &  = \left\{ ( \{\xi,\gamma\}, \{\xi,\beta\} ) ~\vert~  \xi, \gamma, \beta \in \Finf \right\}  
\\\noalign{\medskip}
R_{-1} & = \left\{ ( \{\xi,\gamma\}, \{\alpha,\beta\} ) ~\vert~  \xi, \alpha, \gamma, \beta \in \Finf, 
\mathrm{cr} (\xi, \gamma; \alpha, \beta ) = -1\right\}  
\\\noalign{\medskip}
R_r & = \left\{ ( \{\xi,\gamma\}, \{\alpha,\beta\} ) ~\vert~  \xi, \gamma, \alpha, \beta \in \Finf, 
\mathrm{cr} (\xi, \gamma; \alpha, \beta ) = r ~\mbox{or}~ r^{-1} \right\}, 
r \in \bF_q\setminus \{0, \pm{1}\}.
\end{array} \right\}
\end{equation}
where $\mathrm{cr}(x,y;z,w) = \frac{(x-z)(y-w)}{(x-w)(y-z)}$ is the cross ratio. 
$R_1$ has valency  $2(q-1)$, and $R_{-1} $ has valency $(q-1)/2$, which is half of that of $R_r$.  

\subsection{The case $q$ odd} \label{s:q-odd}
Throughout this subsection, $q$ will be an odd prime power.  The quadratic form $Q(x_0, x_1, x_2) = x_1^2 - x_0 x_2$ polarizes to the symmetric  bilinear form 
$$
B((x_0, x_1, x_2), (y_0,y_1, y_2)) = 2x_1 y_1  - x_0 y_2 - x_2 y_0.
$$  The form $B$ (or the conic $\mathscr{O}$) defines the following \emph{polarity} on $\PG(2,q)$: 
\begin{equation} \label{e:polarity}
\perp:    (x_0: x_1: x_2)  \mapsto  (x_0: x_1, x_2)^{\perp} :=\{ (y_0: y_1: y_2) ~\vert~  B((x_0, x_1, x_2), (y_0,y_1, y_2))  = 0\}.
\end{equation}

Denote by $L_{\xi,\gamma}$  the  hyperbolic line to $\mathscr{O}$ through distinct points $P_\xi$ and $P_\gamma$ for $\xi, \gamma \in \Finf$. It is not difficult to check the following:
\begin{equation} \label{e:perp}
\begin{array}{c}
L_{\xi,\gamma}^\perp = 
\begin{cases} 
 ( \gamma \xi : (\xi + \gamma)/2: 1) & ~\mbox{if}~ \xi, \gamma \in \bF_q,\; \xi \ne \gamma \\
 (2\xi, : 1 : 0)   & ~\mbox{if}~ \xi \in \bF_q, \gamma = \infty,
 \end{cases}
\\
Q(L_{\xi, \gamma}^\perp ) = \left( \frac{\gamma - \xi}{2}\right)^2 \qquad \mbox{and} \qquad 
Q(L_{\xi, \infty}^\perp ) =1.
\end{array}
\end{equation}
A  point $P= (x_0: x_1 : x_2)$  is called {\em hyperbolic} (respectively {\em elliptic})   if $(x_0,x_1,x_2)^\perp$ is is a hyperbolic (respectively elliptic) line. Because of (\ref{e:perp}),  hyperbolic (elliptic) points are also referred to as  \emph{square type} (\emph{non-square type})  in the literature.  Since there are $q(q+1)/2$ secant lines, there are  $q(q+1)/2$ hyperbolic points in total.   Note that $\ell^\perp$ is a elliptic (respectively singular) point if  $\ell$ is an elliptic (respectively tangent) line.  

Denote by $\mathcal{L}_{\epsilon}^{\perp}$ the set of all hyperbolic, singular and elliptic  points for $\epsilon = +, 0,$ and $-$, respectively, Then the action of $\SO(3,q)$  on $\mathcal{L}_{\epsilon}^{\perp}$ determines  a symmetric association scheme, denoted by $\mathfrak{X}(\SO(3,q), \mathcal{L}_{\epsilon}^{\perp})$. We have  the following result.  

\begin{thm} \label{t:scheme-perp}
Let $q$ be odd. Then $\mathfrak{X}(\SO(3,q), \mathcal{L}_{\epsilon})$ is isomorphic to 
$\mathfrak{X}(\SO(3,q), \mathcal{L}_{\epsilon}^{\perp})$. 
\end{thm}

To our best knowledge, $\mathfrak{X}(\SO(3,q), \mathcal{L}_{\epsilon}^{\perp})$ for $\epsilon = \pm$ was first constructed by  Shen \cite{Shen84}. In fact, he studied  association schemes coming from the action of $\SO(n,q)$ on hyperbolic and elliptic points in 
$\PG(n-1,q)$(also cf. \cite{Bannai_Shen_Song90}). Kwok \cite{Kwok91} calculated the character tables of $\mathfrak{X}(\SO(3,q), \mathcal{L}_{\epsilon}^{\perp}) _{\epsilon = \pm}$. 
 We mention in passing that $\mathfrak{X}(\SO(3,q), \mathcal{L}_0)$ is a 2-class symmetric scheme. 

\section{Three fission schemes of $\mathrm{T}(q+1)$} 
\label{s:lattice}

This section we are concerned with the three fission schemes of   $\mathrm{T}(q+1)$ mentioned in Section \ref{s:introduction}.
The next theorem holds for $\PSL(2,q)$ regardless the parity of $q$ (e.g, \cite[Theorem 4.1]{Taylor92}).   

\begin{thm}\label{t:psl2q}
The group $\PSL(2,q)$ acts doubly transitively on $\PG(1,q)$. Hence, $\POm(3,q)$ acts doubly transitively on $\mathscr{O}$.
\end{thm}

Recall that $\Omega$ denotes the collection of 2-element subsets in $\PG(1,q)$. 
By  the above theorem, $\PSL(2,q)$ acts transitively  on $\Omega$, and thus the action  determines an association scheme, denoted by $\mathfrak{X}(\PSL(2,q),   \Omega)$.  This is a fission scheme of $\FT(q+1)$. 

  In the rest of this paper, we shall describe  association schemes from the actions of $\PSL(2,q)$, $\Mq(q)$ and $\PML(2,q)$ on $\Omega$.  In Section \ref{s:isomorphism}, we give a  one-to-one correspondence  of $\Omega$  with $\mathcal{L}_+$ (and   
  $\mathcal{L}_+^\perp$ via polarity)  and  an embedding $\rho$ of $\PML(2,q)$.  As result, for any transitive subgroup of $\PML(2,q)$ on $\Omega$, 
  we obtain three isomorphic association schemes, e.g., $\mathfrak{X}(\PSL(2,q),   \Omega)$,  $\mathfrak{X}(\POm(3,q), \mathcal{L}_+)$ and $\mathfrak{X}(\POm(3,q), \mathcal{L}_+^\perp)$ in the case . of $\PSL(2,q)$.
So we shall use the action on  $ \mathcal{L}_+^\perp$ in subsequent calculation.  

\subsection{The association scheme $\mathfrak{X}(\POm(3,q), \mathcal{L}_+)$}  \label{s:po3q}
In this subsection, we consider the association scheme from the action of $\POm(3,q)$. Namely, we shall prove the following result. 

\begin{thm}\label{t:pom}
The association scheme $\mathfrak{X}(\POm(3,q), \mathcal{L}_{+})$ is non-symmetric with $(3q+5)/4$ classes if
 $q \equiv 1 \bmod 4$ and $(3q+3)/4$ classes if  $q \equiv 3 \bmod 4$.
\end{thm}

Consider the stabilizer of  $L_{0,\infty}^\perp$ in $\POm(3,q)$, which is denoted by 
$\POm(3,q)_{0,\infty}$.
We know from Subsection \ref{s:embedding} that $\POm(3,q)$ is induced by elements of form (\ref{e:so3q}) with $ad-bc=1$.  Since  $L_{\infty, 0}^\perp$ has coordinates $(0:1:0)$ by (\ref{e:perp}), 
$$
\begin{pmatrix}  a^2 & 2ab & b^2 \\  ac & ad +bc & bd \\ c^2 & 2cd & d^2 \end{pmatrix}  
\col{ 0\\ 1 \\ 0} =   \col{2ab \\ ad+bc  \\ 2cd } = \col{0 \\ \lambda \\ 0}
$$
for some $\lambda$ in $\bF_q^*$.  Since $ad - bc = 1$, we must have  $a= d^{-1},  b =c =0$, or $a=d=0, b = -c^{-1}$. So $G := \POm(3,q)_{0,\infty}$ is induced by the following elements 
\begin{equation}
\begin{pmatrix}  a^2 & & \\ & 1 & \\ &  & a^{-2} \end{pmatrix},  \quad
\begin{pmatrix} & & b^{-2}  \\ & -1 & \\ b^2 & & \end{pmatrix} ,\quad  a, b \in \bF_q^*. 
\end{equation}

For each relation $R$ in (\ref{e:rels}), let 
$$
R^{\{0,\infty\}} = \left\{ L_{\xi, \gamma} ~\vert~  (L_{0,\infty}, L_{\xi, \gamma}) \in R
\right\}.
$$
These $R^{\{0,\infty\}} $ are all orbits of $\SO(3,q)$ on $\mathcal{L}_{+}$.  
Since  $\POm(3,q)$ is a subgroup of $\SO(3,q)$ of index 2.  So any orbit of $\SO(3,q)$   either splits into a pair of orbits  or remains to be an orbit of $\POm(3,q)$.    

{\em
In the rest of this paper, $z$ is a fixed non-square element in $\bF_q$, $s$ denotes a square element  not in   $\{0, \pm{1}\}$, and $t$ denotes a non-square element in $\bF_q$. Let $\bF_q^{*2}$ be the set of all nonzero square elements in $\bF_q$. So $z\bF_q^{*2}$ consists of all non-square elements of $\bF_q$.
}

By (\ref{e:perp}), the hyperbolic points come in with two forms:
\begin{equation}
(2\xi : 1 : 0) ~\mbox{with $\xi \in \bF_q$}, \quad ~\mbox{or}~ \quad 
(\gamma\xi: (\xi + \gamma)/2: 1) 
~\mbox{with $\xi,\gamma \in \bF_q,\ \xi\ne \gamma$}.
\end{equation}

First, we settle the $G$-orbits on $(2\xi : 1 : 0)$.  
It is trivial that    $\Gamma_0 =  \{ L_{0, \infty}\} $ is a $G$-orbit. Now assume $\xi \ne 0$.  
Note that
\begin{equation}
\begin{pmatrix}  a^2 & & \\ & 1 & \\ &  & a^{-2} \end{pmatrix} \col{2\xi \\  1 \\ 0 } = \col{2a^2\xi \\ 1 \\ 0 },  \quad
\begin{pmatrix} & & a^{-2}  \\ & -1 & \\ a^2 & & \end{pmatrix}  \col{2\xi \\  1 \\ 0 } = 
\col{0 \\ -1 \\ 2a^2 \xi } = \col{ 0 \\ \frac{-1}{2a^2\xi} \\ 1}, 
\end{equation}
where the last equality follows from that these two vectors represent the same (projective) point. 
Recall that  hyperbolic points are in one-to-one correspondence with 2-element subsets in $\mathscr{O}$, which are again  in one-to-correspondence with those of $\PG(1,q)$.   Note that 
$(2a^2\xi : 1 : 0 )$  is determined by  $\{a^2\xi, \infty\}$, and $( 0 : \frac{-1}{2a^2\xi} : 1)$  by $\{0,-(a^2\xi)^{-1}\}$.
For a fixed square  (respectively non-square) $\xi$  in $\bF_q^*$,  if $-1\in {\bF_q^*}^2$,  then both $a^2\xi$ and  $-(a^2\xi)^{-1}$ range over ${\bF_q^*}^2$ (respectively  $z\bF_q^{*2}$).  On the other hand,  if $-1\notin {\bF_q^*}^2$,  then $a^2\xi$ ranges over ${\bF_q^*}^2$  while $-(a^2\xi)^{-1}$ ranges over   $z\bF_q^{*2}$.
Hence,  $R_1^{\{0,\infty\}}$ splits into two orbits of   length $q-1$ each:
\begin{equation}\label{e:orb1}
\left\{
\begin{array}{lll}
\Gamma_1^{+} =\left\{L_{0, \xi} , L_{\infty, \xi}  ~\vert~ \xi  \in  {\bF_q^*}^2 \right\}  &
\Gamma_1^{-} = \left \{L_{0, \xi} , L_{\infty, \xi} ~\vert~ \xi  \in  z\bF_q^{*2} \right \}
& ~\mbox{ if $q\equiv 1 \bmod 4$} 
\\ [5pt]
\Gamma_1^{+} =\left\{L_{0, \xi} , L_{\infty, \gamma}  ~\vert~ \xi  \in  {\bF_q^*}^2, \, \gamma \in z{\bF_q^*}^2 \right\} 
&
\Gamma_1^{-} =\left\{L_{0, \xi} , L_{\infty, \gamma}  ~\vert~ \xi  \in  {z\bF_q^*}^2, \, \gamma \in {\bF_q^*}^2 \right\}  
& ~\mbox{  if $q\equiv 3 \bmod 4$}.
\end{array}\right.
\end{equation}

The remaining  hyperbolic points are $(\gamma\xi: (\xi +\gamma)/2: 1)$ with $\xi \gamma \ne 0.$ 
Note that  
\begin{equation} \label{e:diag}
\begin{pmatrix}  a^2 & & \\ & 1 & \\ &  & a^{-2} \end{pmatrix}
 \col{
 \gamma\xi 
  \\\noalign{\smallskip} \frac{\xi +\gamma}{2} 
  \\\noalign{\smallskip} 1} 
 = 
 \col{
 \gamma_1\xi_1 
  \\\noalign{\smallskip} \frac{\xi_1 +\gamma_1}{2} 
  \\\noalign{\smallskip} 1}  
  \quad \mbox{if and only if } \quad 
\{  a^2 \gamma,   a^2 \xi\} =\{ \gamma_1,  \xi_1\}
\end{equation}
and 
\begin{equation} \label{e:offdiag}
\begin{pmatrix} & & b^{-2}  \\ & -1 & \\ b^2 & & \end{pmatrix}  
\col{\gamma\xi
 \\\noalign{\smallskip} \frac{\xi +\gamma}{2} 
  \\\noalign{\smallskip} 1} 
  = 
  \col{
 \gamma_1\xi_1 
  \\\noalign{\smallskip} \frac{\xi_1 +\gamma_1}{2} 
  \\\noalign{\smallskip} 1}  
  \quad \mbox{if and only if } \quad 
  \left\{ \dfrac{-1} { b^2 \gamma},  \dfrac{-1}{b^2\xi} \right\} =\{    \gamma_1,   \xi_1  \}.
\end{equation}

If $\gamma/\xi = -1$, $R_{-1}^{\{0,\infty\}}$ remains to be a $G$-orbit for $q\equiv 3\bmod 4$ and it splits into two orbits for $q\equiv 1\bmod 4$:
\begin{equation}  \label{e:orb2}
 \left\{
\begin{array}{cl}
\Gamma_{-1} =  \left\{L_{\xi, -\xi} ~\vert~ \xi \in {\bF_q^*}\right\} & \mbox{if $q\equiv 3\bmod 4$}
 \\\noalign{\smallskip}
\Gamma_{-1}^{+} = \left\{L_{\xi, -\xi} ~\vert~ \xi \in {\bF_q^*}^2\right\},
\quad
\Gamma_{-1}^{-} = \left\{L_{\xi, -\xi} ~\vert~ \xi \in \bF_q^* \setminus {\bF_q^*}^2\right\} &  \mbox{if $q\equiv 1\bmod 4$}.
\end{array}\right.
\end{equation}
The orbit length is  $(q-1)/2$  in the first case and $(q-1)/4$ in the second case. 

Now suppose $r = \gamma/\xi \notin \{0, \pm{1}\}$. We claim that $R_\gamma^{\{0, \infty\}}$  becomes the following $G$-orbit(s):
\begin{equation} \label{e:orb3}
\left\{
\begin{array}{cl}
\Gamma_{r} =  \left\{L_{\xi, r\xi}~\vert~ \xi \in {\bF_q^*}\right\} & \mbox{if $-1/r$ is non-square}
 \\\noalign{\smallskip}
\Gamma_{r}^{+} = \left\{L_{\xi, r\xi} ~\vert~ \xi \in {\bF_q^*}^2\right\},  
\quad
\Gamma_{r}^{-} =  \left\{L_{\xi, r\xi}~\vert~ \xi \in {z\bF_q^*}^2\right\}&  \mbox{if $-1/r$ is square}.
\end{array}\right.
\end{equation}
The orbit length is $(q-1)$ in the first case and  $(q-1)/2$ in the second case.  

Now suppose that $\gamma= r\xi$ and $\gamma_1 = r\xi_1$. If $\xi$ and $\xi_1$ are both squares (non-squares), set 
$a^2 = \xi_1\xi^{-1}=a^2$. By (\ref{e:diag}),  $L_{\xi, r\xi}$ and $L_{\xi_1, r\xi_1}$ are in the same $G$-orbit.
  Suppose that  $\xi$ is a square and $\xi_1$ is 
a non-square. If $-1/r$ is a non-square, set $b^2 = \frac{-1}{r\xi\xi_1}$. Then   $L_{\xi, r\xi}$ and $L_{\xi_1, r\xi_1}$ are in the same
 $G$-orbit. 
 On the other hand,  if $-1/r$ is a square, then $\{\xi_1, r\xi_1\} = \left\{\frac{-1}{b^2\xi}, \frac{-1}{b^2 r\xi }\right\}$ if and
   only if $r = \pm {1}$. As $r\ne \pm{1},$   $L_{\xi, r\xi}$ and $L_{\xi_1, r\xi_1}$ fall into different orbits.  We prove the above claim. 

Each orbit $\Gamma$ determines a  relation $R$ of $\mathfrak{X}(\POm(3,q), \mathcal{L}_{+})$, and vice versa.   So we  may equip  $R$  with the same subscript/superscript of $\Gamma$, e.g., $R_{-1}$, $R_{ r}^+$, etc.  
We will adopt this convention for the other two fission schemes to follow. 

Now we determine the non-symmetric relations of  $\mathfrak{X}(\POm(3,q), \mathcal{L}_{+})$.  As an example, we illustrate with 
${R}_{s}^{+}$.    We shall use $\PSL(2,q)$, instead of $\POm(3,q)$. Take $(L_{0, \infty}, L_{1, s}) \in R_s^+$. Note the function $x \mapsto \frac{x-1}{x-s}$ maps $\{1, s\}$ to  $\left\{ 0, \infty\right\}$.
Now suppose that  $\det \begin{pmatrix} 1 & -1 \\ 1 & -s\end{pmatrix} = 1-s$ is a non-square, say $z$. The function
 $x \mapsto \frac{x -1}{z(x -s)} $ is in $\PSL(2,q)$. This function carries $\{1,s\}$ to $\{0,\infty\}$, and $\{0, \infty\}$ to 
 $\{\frac{1}{z}, \frac{1}{sz}\}$.  The  hyperbolic line corresponding to  $\{\frac{1}{z}, \frac{1}{sz}\}$ is in $\Gamma_s^-$ and thus 
 ${^tR}_{s}^{+}  =     R_{s}^{-1}$. 
We can handle other relations in a similar way.  

 Suppose that $q \equiv 1 \bmod 4$. Then ${^tR}_{-1}^{+}  =       R_{-1}^{-}$  if $2$ is a non-square in $\bF_q$, 
and ${^tR}_{s}^{+}  =     R_{s}^{-1}$   if  $1-s $ is a non-square.
It is well known that $2$ is a square in $\bF_q$  if and only if the Legendre symbol 
$\left(\frac{2}{q}\right)  = \left[(-1)^{\frac{p^2-1}{8}} \right]^m =1$, where $q = p^m$ for some  odd prime $p$ (e.g., see 
\cite[Section 6.1]{Ireland}). 
Suppose that $q \equiv 3 \bmod 4$. Then 
$
{^tR}_1^{+} = R_1^{-}$, 
and 
${^tR}_{t}^{+}  =   R_{t}^{-1}$  if  $1-t$ is a square.

By counting orbits in (\ref{e:orb1}), (\ref{e:orb2}) and (\ref{e:orb3}), we see that    $\mathfrak{X}(\PML(q), \mathcal{L}_{+})$ has the  asserted class number. The proof of Theorem \ref{t:pom} is complete.

 \subsection{The association scheme $\mathfrak{X}(\Mq(q), \mathcal{L}_{+})$}

 In this subsection, we describe the association scheme  $\mathfrak{X}(\Mq(q), \mathcal{L}_{+})$, which is a fusion scheme of 
 $\mathfrak{X}(\POm(3,q), \mathcal{L}_{+})$.  Now we prove the following result.

\begin{thm} \label{t:mq}
Let $ q$ be an even power of an odd prime. The following holds:
\begin{enumerate}
\item \label{t:mqg}
For  $q>9$,  $\mathfrak{X}(\Mq(q), \mathcal{L}_{+})$ is a non-symmetric association scheme with  $(3q+5)/8$ classes ;
\item  \label{t:mq9}
$\mathfrak{X}(\Mq(9), \mathcal{L}_{+})$ is a symmetric $P$-polynomial association scheme.  
\end{enumerate}
 \end{thm}

 Now consider the stabilizer $\Mq(q)_{0,\infty}$ of $L_{0,\infty}^\perp$ in $\Mq(q)$. 
 Under the embedding $\rho$ of $\PML(2,q)$ into $\PML(3,q)$ in Subsection \ref{s:embedding},   $\rho(S) = \POm(3,q)$.  The stabilizer of $ L_{0,\infty}$ in $\Mq(q)$ is the union of $\POm(3,q)_{0,\infty}$ and $\rho(T)_{0,\infty}$, where 
 $\rho(T)_{0,\infty}$ consists of elements in $\rho(T)$ that fix $L_{0,\infty}$. Then  $\rho(T)_{0,\infty}$ is induced by all the   semilinear transformations of the following matrices and  the unique involution $\sigma$:
 \begin{equation}
 \begin{pmatrix} a^2 & \\ & \zz& \\ && a^{-2}\zz^{2} \end{pmatrix}, \quad
 \begin{pmatrix}  && b^2 \\ & -\zz& \\ b^{-2}\zz^2&& \end{pmatrix} ,
 \end{equation}
 where $a, b$ ranges over ${\bF_q^*}$. 
 
For typographical convenience, we write $\tilde{x} = x^\sigma$. Since $\tilde{x} = x$ for any $x$ in the fixed field 
$\mathrm{Fix}(\sigma)$ of $\sigma$, we can check that,  for any points $(2\theta: 1:0)$ and $(\gamma\theta : (\theta+\gamma)/2:1)$, 
  \begin{equation} \label{e:Mq1}
   \begin{pmatrix} a^2 & \\ & \zz& \\ && a^{-2}\zz^{2} \end{pmatrix}
   \begin{pmatrix} 2\tilde{\theta} \\  1  \\  0  \end{pmatrix} 
   =
     \begin{pmatrix}  2a^2  \tilde{\theta} \zz^{-1} \\ 1 \\ 0 \end{pmatrix}, 
     \quad
   \begin{pmatrix}  && b^2 \\ & -\zz& \\ b^{-2}\zz^2&& \end{pmatrix} 
     \begin{pmatrix} 2\tilde{\theta} \\  1  \\  0  \end{pmatrix} 
 =   \begin{pmatrix}  0 \\ \frac{-b^2}{2\zz  \tilde{\theta} } \\ 1  \end{pmatrix}, 
  \end{equation}
 We remark that $\mathrm{Fix}(\sigma)$ has $\sqrt{q}$ elements. 
  
 Note a transformation in (\ref{e:Mq1}) maps $\Gamma_1^{+}$ to $\Gamma_1^{-}$, and vice versa. So 
 $\Mq(q)_{0, \infty}$ combines $\Gamma_1^{+}$ and $ \Gamma_1^{-}$:
 $$
 \Delta_1 = \Gamma_1^{+} \cup \Gamma_1^{-}.
 $$
 The following equations shall be used in determining the remaining orbits of  $\Mq(q)_{0, \infty}$:
  \begin{equation} \label{e:Mq2}
      \begin{pmatrix} a^2 & \\ & \zz& \\ && a^{-2}\zz^{2} \end{pmatrix}
   \begin{pmatrix} \tilde{\gamma} \tilde{\theta} \\ \frac{ (\tilde{\theta}+ \tilde{\gamma})}{2} \\ 1 \end{pmatrix} 
    =  \begin{pmatrix} \frac{a^4\tilde{\gamma} \tilde{\theta}} {\zz^{2}} \\[5pt] \frac{a^2(\tilde{\theta}+ \tilde{\gamma})}{2\zz} \\ 1 \end{pmatrix},  
\quad 
    \begin{pmatrix}  && b^2 \\ & -\zz& \\ b^{-2}\zz^2&& \end{pmatrix} 
         \begin{pmatrix} \tilde{\gamma} \tilde{\theta} \\ \frac{ (\tilde{\theta}+ \tilde{\gamma})}{2} \\ 1 \end{pmatrix} 
     = 
     \begin{pmatrix}  \frac  {b^4} { \tilde{\theta} \tilde{\gamma} \zz^2} \\[5pt]
      -\frac{(\tilde{\theta} + \tilde{\gamma})b^2}{2 \tilde{\theta}\tilde{\gamma}\zz}  \\ 1\end{pmatrix}.
  \end{equation}
 We can deduce from (\ref{e:Mq2})  that  $\Mq(q)_{0, \infty}$ combines $\Gamma_{-1}^+$ and $\Gamma_{-1}^{-}$: 
  $$\Delta_{-1} =  \Gamma_{-1}^+\cup \Gamma_{-1}^{-}.
  $$
 Note that  $\Delta_1$ and $\Delta_{-1}$  have lengths $2(q-1)$ and $(q-1)/2$, respectively. 
  
 Now we consider the orbits of   $\Mq(q)_{0, \infty}$ on the remaining points $(\theta\gamma: (\theta+\gamma)/2: 1)$ with 
 $\gamma + \theta \ne 0$. Let $r = \theta / \gamma.$  These orbits depend on whether $\tilde{r}$ is in $\{r, r^{-1}\}$ or not.  Note if $\tilde{r} = r $ (respectively $r^{-1} $) , then $r^{\sqrt{q} -1} = 1$ (respectively   $r^{\sqrt{q} +1} = 1$) and thus $r$ is a square in $\bF_q$. In this case, we obtain $(\sqrt{q}-3)/2$ (respectively $(\sqrt{q}-1)/2$) orbits of length $(q-1)$ each. Now change $r$ to $s$. 
 So we can deduce  from (\ref{e:Mq2}) that   $\Mq(q)_{0, \infty}$ combines $\Gamma_{s}^+$ and 
 $\Gamma_{s}^-$ if $ \tilde{s} = s$ or $s^{-1}$:
 $$
\Delta_{s}  =  \Gamma_{s}^+ \cup   \Gamma_{s}^-.
 $$

 If  $ s\in \bF_q^{*2}$ and $\tilde{s} \notin \{s, s^{-1}\}$, then we can deduce again from (\ref{e:Mq2})  that  $\Mq(q)_{0, \infty}$ combines $\Gamma_{s}^+$ and 
 $\Gamma_{\tilde{s}}^-$ (respectively  $\Gamma_{s}^-$ and 
 $\Gamma_{\tilde{s}}^+$):
 $$
 \Delta_s^+ = \Gamma_{s}^+ \cup  \Gamma_{\tilde{s}}^-, 
 \qquad 
  \Delta_s^- =   \Gamma_{s}^- \cup \Gamma_{\tilde{s}}^{+} .  
 $$ 
 We obtain $\frac{(\sqrt{q}-3)(\sqrt{q}-1)}{4}$ orbits of length $q-1$. 
 
If $t$ is a non-square in $\bF_q$, then  $\Mq(q)_{0, \infty}$ combines $\Gamma_{t}$ and $\Gamma_{\tilde{t}}$:
$$
\Delta _{t} = \Gamma_{t} \cup \Gamma_{\tilde{t}}. 
$$
In this case, we obtain $(q-1)/8$ orbits of length $2(q-1)$ each. 

Now we write the relations of $\mathfrak{X}(\Mq(q), \mathcal{L}_{+})$ according to these $\Delta$.   If $q\ne 9$,
$\mathfrak{X}(\Mq(q), \mathcal{L}_{+})$ is non-symmetric:   ${^tR}_s^+ = R_s^-$.  So we prove (\ref{t:mqg}) of Theorem \ref{t:mq}.
   $\mathfrak{X}(\Mq(9), \mathcal{L}_{+})$ is a symmetric scheme; see Example \ref{e:exam} below. In fact, it is a $P$-polynomial scheme  on 45 vertices  from a generalized octagon of order $(2,1)$ \cite[p. 419]{BCN}.  We complete the proof of Theorem 
    \ref{t:mq}.

\begin{rmk}  Limited computation with Gap \cite{Gap} shows that $\mathfrak{X}(\Mq(q), \mathcal{L}_{+})$ is not commutative for $q>25$. 
it is not symmetric but commutative for $q=25$. 
  
\end{rmk}

\subsection{The association scheme $\mathfrak{X}(\PML(2,q), \mathcal{L}_{+})$}  

In this subsection, we consider an association scheme that is a fusion scheme of both $\FT(q+1)$ and
$\mathfrak{X}(\Mq(q), \Omega)$.  Since both $\PGL(2,q)$ and $\Mq(q)$ are subgroups of 
$\PML(2,q)$, the action of $\PML(2,q)$ on $\Omega$ determines this fusion scheme. 
 
The stabilizer of $\rho(\PML(2,q))_{0, \infty}$ of $L_{0,\infty}^\perp$ induced by all the  semilinear  
 transformation of the following matrices associated with an automorphism $\tau$ of $\bF_q$:
  \begin{equation}
 \begin{pmatrix} a^2 & \\ & ad & \\ && d^{2} \end{pmatrix}, \quad
 \begin{pmatrix}  && b^{2} \\ & bc & \\ c^2 && \end{pmatrix} ,
 \end{equation}
 where $a, b,c,d$ range over ${\bF_q^*}$, and $\tau$ ranges over $\mathrm{Aut}(\bF_q)$ . 
 
The group $\rho(\PML(2,q))_{0, \infty}$ acts on $\mathcal{L}_{+}^\perp$ similarly to that in (\ref{e:Mq1}) and  (\ref{e:Mq2}) with $\;\tilde{}\;$ replaced with $\tau$.  By a  similar   argument employed in the previous subsection, $\rho(\PML(2,q))_{0, \infty}$ has the following orbits on $\mathcal{L}_{+}$:
\begin{equation} \label{e:Lam}
\begin{gathered}
\Lambda_{0} = \{L_{0, \infty}\}, \quad \Lambda_{1} = \{L_{0,\xi}, L_{\infty, \xi} ~\vert~ \xi \in \bF_q^*\}, 
\quad  \Lambda_{-1} = \{ L_{\xi, -\xi} ~\vert~ \xi \in \bF_q^*\}, \\
\Lambda_{r} = \{L_{\xi, \gamma} ~\vert~ \xi / \gamma\in \{r^\tau, r^{-\tau} ~\vert~ \tau \in \mathrm{Aut\;}\bF_q\}\}, (r\in \bF_q^*\setminus \{ \pm{1}\}).
\end{gathered}
\end{equation}
No we can write down the corresponding relations of  $\mathfrak{X}(\PML(2,q), \mathcal{L}_{+})$:
\begin{equation*}
\begin{gathered}
R_0 = \{(L_{\xi, \gamma}, L_{\xi, \gamma}) ~\vert~ \{\xi,  \gamma\} \subset \Finf \},
\\
R_1 =\{ (L_{\xi, \gamma}, L_{\xi, \beta}) ~\vert~ \{\xi, \gamma, \beta\} \subset \Finf \}, \\
R_{-1} =   \{(L_{\xi, \gamma}, L_{\alpha, \beta}) ~\vert~  \{\xi, \gamma, \alpha, \beta\} \subset \bF_q^*,\ 
\mathrm{cr}(\xi, \gamma; \alpha, \beta) = -1\},
\\
R_r = \left\{(L_{\xi, \gamma}, L_{\alpha, \beta}) ~\vert~  \{\xi, \gamma, \alpha, \beta\} \subset \bF_q^*,\ 
\mathrm{cr}(\xi, \gamma; \alpha, \beta) \in  \{r^{\tau}, r^{-\tau} ~\vert~ \tau \in \mathrm{Aut\;}\bF_q\}\right\},
 \qquad r\in \bF_q^*\setminus \{ \pm{1}\}.
\end{gathered}
\end{equation*}

\begin{thm}
 $\mathfrak{X}(\PML(2,q), \mathcal{L}_{+})$ is a symmetric association scheme with the above relations. 
\end{thm}

 \begin{rmk} 
  Note that the relation $R_1$ is an original relation in $\mathrm{T}(q+1)$.
For $q=9$,  $\mathfrak{X}(\PML(2,q), \mathcal{L}_{+})$ has 4 classes.  It has the $P$-polynomial property given by $R_{-1}$ but not the $Q$-polynomial property.  For $q =25$ and $ 49$,  $\mathfrak{X}(\PML(2,q), \mathcal{L}_{+})$ has $9$ and $16$ classes, respectively.  These schemes don't have   polynomial ($P$ or $Q$) property.  For the general case, the class number  of   $\mathfrak{X}(\PML(2,q), \mathcal{L}_{+})$ is reduced to   counting  the number of orbits of $\mathrm{Aut}(\bF_q)$ on the collection 
$\{ \{r^\tau, r^{-\tau}\} ~\vert~  r\in \bF_q^*\setminus \{ \pm{1}\}, \; \tau \in \mathrm{Aut}(\bF_q)\ \}.$ 
However, we are unable to give  an explicit formula for this number.  
\end{rmk}
 
\begin{example} \label{e:exam} {\em
In this example, we illustrate the four fission schemes of $\mathrm{T}(10)$ for $q =9$. 
Let $g$ be a primitive element of $\bF_9$. We have  
$\mathfrak{X}(\PML(2,9), \mathcal{L}_{+})= \mathfrak{X}(\Mq(9),\mathcal{L}_{+})$. 
The relations of $\FT(10), \mathfrak{X}(\PSL(2,9), \mathcal{L}_{+})$ and  
$\mathfrak{X}(\Mq(9), \mathcal{L}_{+})$ are shown in the following diagram:

$$
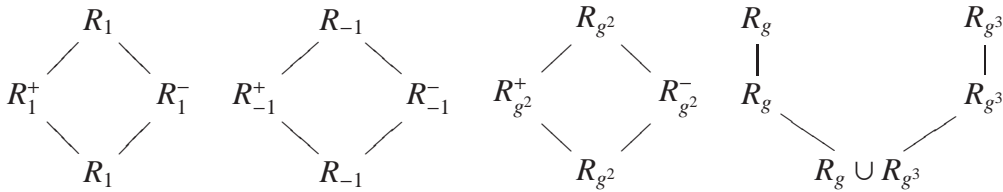
\xymatrix@-1.2pc{
  & R_1 \ar@{-}[dl] \ar@{-}[dr] & & & R_{-1}  \ar@{-}[dl] \ar@{-}[dr]& & &R_{g^2}  \ar@{-}[dl] \ar@{-}[dr]& 
& R_g\ar@{-}[d] & & R_{g^3} \ar@{-}[d]  \\
 R_1^+ && R_1^- & R_{-1}^+ && R_{-1}^- & R_{g^2}^+ && R_{g^2}^- & R_g &&R_{g^3}    \\
  & R_1 \ar@{-}[ul] \ar@{-}[ur]  & & & R_{-1}\ar@{-}[ul] \ar@{-}[ur]  & & &R_{g^2} \ar@{-}[ul] \ar@{-}[ur] && & R_g \cup  R_{g^3}\ar@{-}[ul] \ar@{-}[ur]  &   \\
}
$$
The top, middle and bottom rows are the relations of $\FT(10), \mathfrak{X}(\PSL(2,9), \mathcal{L}_{+})$ and  
$\mathfrak{X}(\Mq(9), \mathcal{L}_{+})$, respectively.  
} \end{example}
 
\begin{rmk}
Primitive representations of $\PSL(2,q)$ were investigated by Sasha Ivanov and Loius Tchuda via the use of so-called 
Burnside marks  \cite{Klin}. 
A general method was outlined in \cite{Klin94}, which allows one to obtain formulas for the rank and subdegrees for a (transitive) permutation representation \cite{Klin}.
\end{rmk}
\section*{Acknowledgment} 
When this project was started, JM had helpful discussions with Professor Robert  Liebler several months before his unexpected death in July 2009.  The authors are grateful to Professor Tim Penttila   for  helpful conversations on numeral occasions, and to Professors Edwin van Dam and   Hajimie Tanaka for their comments on  a draft.   
 K. Wang's research is supported by NCET-08-0052, NSF of China (10871027) and the Fundamental Research Funds for
 the Central University of China.


\end{document}